\documentclass[9pt]{jdg-p-1-2}%
\usepackage{amsmath}%
\usepackage{amsfonts}%
\usepackage{amssymb}%
\usepackage{graphicx}
\providecommand{\U}[1]{\protect\rule{.1in}{.1in}}
\newtheorem{theorem}{Theorem}[section]
\newtheorem{lemma}[theorem]{Lemma}
\newtheorem{corollary}[theorem]{Corollary}
\newtheorem{proposition}[theorem]{Proposition}
\theoremstyle{definition}

\theoremstyle{remark}

\def\Ric{\text{Ric}}

\def\e{\epsilon}
\def\p{\partial}

\def\R{\Bbb R}

\def\id{\operatorname{id}}
\def\Ric{\operatorname{Ric}}
\def\I{\operatorname{I}}
\def\Scal{\operatorname{Scal}}
\def\tr{\operatorname{tr}}

\def\W{\operatorname{W}}
\def\l{\operatorname{l}}

\def\tl{\tilde{\lambda}}

\numberwithin{equation}{section}
\begin{document}

\title{Complete manifolds with nonnegative curvature operator }
\author{Lei Ni}\thanks{The first author was supported in part by NSF Grants and an Alfred P. Sloan
Fellowship, USA}



\address{Department of Mathematics, University of California at San Diego, La Jolla, CA 92093}

\email{lni@math.ucsd.edu}



\author{Baoqiang Wu}

\address{Department of Mathematics, Xuzhou Normal University, Xuzhou, Jiangsu, China}

\email{wubaoqiang@xznu.edu.cn}

\date{June 2006}

\begin{abstract}

In this short note, as a simple application of the strong result
proved recently by B\"ohm and  Wilking, we give a classification on
closed manifolds with $2$-nonnegative curvature operator. Moreover,
by the new invariant cone constructions  of B\"ohm and  Wilking, we
show that any complete Riemannian manifold (with dimension $\ge 3$)
whose curvature operator is bounded and satisfies the pinching
condition $R\ge \delta R_{I}>0$, for some $\delta>0$, must be
compact. This provides an intrinsic analogue of a result of Hamilton
on convex hypersurfaces.
\end{abstract}

\maketitle

\section{Introduction}
 Let $(M, g)$ be a Riemannian manifold. The curvature operator of $(M, g)$
 lies in the subspace $S_B^2(\wedge^2TM)$ of
  $S^2(\wedge^2TM)$ cut out  by the Bianchi identity. The
  decomposition $S_B^2(\wedge^2TM)=\langle \I\rangle \oplus \langle
  \Ric_0\rangle \oplus \langle \W\rangle$ splits the space of algebraic curvature
  operators into $O(n)$-invariant orthogonal irreducible subspaces. For an orthonormal
basis $\phi_{\alpha}$ (say $\phi_\alpha=e_i\wedge e_j$) of
$\wedge^2TM$ (which can be identified with $so(n)$), the Lie bracket
is given in terms of
$$
[\phi_\alpha, \phi_\beta]=c_{\alpha \beta \gamma} \phi_\gamma.
$$
It is easy to check, by simple linear algebra, that $ \langle [\phi,
\psi],\omega\rangle =-\langle [\omega, \psi], \phi\rangle. $ Here
$\langle A, B\rangle=-\frac{1}{2}\tr(AB)$. This immediately implies
that $c_{\alpha \beta\gamma}$ is anti-symmetric. If $A, B\in
S^2(\wedge^2TM)$ one can define
$$
(A\#B)_{\alpha \beta}=\frac{1}{2}c_{\alpha \gamma \eta}c_{\beta
\delta \theta}A_{\gamma \delta}B_{\eta \theta}.
$$
It is easy to see that $A\#B$ is symmetric too. Also from the
anti-symmetry of $c_{\alpha\beta\gamma}$
$
A\#B=B\#A.
$

 In \cite{BW}, a remarkable algebraic identity was proved on how a linear
 transformation of $S_B^2(\wedge^2 TM)$ changes the quadratic
 form $Q(R)=R^2+R^\#$.
 B\"ohm and  Wilking then constructed a continuous {\it pinching family of invariant closed convex
 cones}. Using this construction they
 confirmed a conjecture of Hamilton stating that {\it
 on a compact manifold the normalized Ricci flow evolves a
 Riemannian metric with $2$-positive curvature operator to a limit
 metric with constant sectional curvature}. Hence it gives a
 complete topological classification of compact manifolds with positive
 $2$-positive curvature operator. In this short notes, based on the
 strong result and the techniques of \cite{BW}, we give the
 classification for manifolds with $2$-nonnegative curvature
 operators and an application of their invariant cone constructions  to
 the compactness of Riemannian manifolds with pinching curvature
 operator.

\section{A strong maximum principle}
Let $(M, g(t))$ be  a complete solution to Ricci flow such that
there exists a constant $A$ and the  curvature tensor of $g(t)$
satisfies $|R_{ijkl}|^2(x,t) \le A$, for all $(x, t)\in M\times[0,
T]$. In \cite{H86}, Hamilton proved that under the evolving normal
frame the curvature tensor satisfies the following evolution
equation.
\begin{proposition}[Hamilton]
\begin{equation}
\label{ham861}\left( \frac{\partial}{\partial t}-\Delta\right)
R=2\left( R^{2}+R^{\#}\right)
\end{equation}
where $R^{\#}=R\#R$.
\end{proposition}

The following was observed for compact manifolds in \cite{Chen,
H90}. We spell out the argument for the noncompact case for the sake
of the completeness.
\begin{proposition}\label{chen1}
The convex cone of $2$-nonnegative curvature operator is preserved
under the Ricci flow.
\end{proposition}
\begin{proof} Let $\I$ be the identity of $S_B^2(\wedge^2TM)$, which
can be identified with the induced metric on $\wedge^2TM$ (as a
section of $\wedge^2TM\otimes \wedge^2TM$). We also denote the
identity map of $TM$ by $\id$. With respect to the evolving normal
frame we have that $\nabla \I=0$ and $\frac{\p}{\p t}\I=0$. Let
$\psi(x,t)> 0$ be the fast growth function constructed in Lemma 1.1
of \cite{NT1} satisfying $\frac{\p}{\p t}\psi-\Delta \psi \ge C_1
\psi$. Here $C_1$ can be chosen as arbitrarily large as we wish. We
shall consider $\tilde R=R+\epsilon \psi \I$ and show that $\tilde R
$ is $2$-positive  for every (sufficiently small) $\epsilon$. If not
by the boundedness of $R$ and growth of $\psi$ we know that it can
only fail somewhere finite. Assume that $t_0$ is the first time
$\tilde R$ fails to be $2$-positive and it happens at some point
$x_0$. If we choose orthonormal basis $\omega_\alpha$ (it may not be
in the form of $e_i\wedge e_j$ as $\phi_\alpha$) such that $\tilde
R$ is diagonal (so is $R$) with eigenvalue $\mu_1\le
\mu_2\le\cdot\cdot \cdot \le \mu_N$, where $N=\frac{n(n-1)}{2}$.
Parallel translate $\omega_\alpha$ to a neighborhood of $(x_0,
t_0)$, and let $\tilde R_{\alpha\alpha}=\langle R(\omega_\alpha),
\omega_\alpha\rangle$ then
 at $(x_0, t_0)$ we have, by the maximum principle,  that
\begin{eqnarray*}
0&\ge& \left(\frac{\p}{\p t} -\Delta \right)\left(\tilde
R_{11}+\tilde R_{22}\right)\\
&\ge& (R^2+R^\#)_{11}+(R^2+R^\#)_{22}+2\epsilon C_1\psi\\
&=& (\tilde R^2+\tilde R^\#)_{11}+(\tilde R^2+\tilde
R^\#)_{22}+2\epsilon C_1\psi\\
&\quad &+\left(R^2+R^\#-\tilde R^2-\tilde
R^\#\right)_{11}+\left(R^2+R^\#-\tilde R^2-\tilde R^\#\right)_{22}
\\
&=&\mu_1^2 +\mu_2^2+\sum
(c^2_{1\beta\gamma}+c^2_{2\beta\gamma})\mu_\beta\mu_\gamma+2\epsilon
C_1\psi
\\&\quad &
-\epsilon \psi\left(\left(2\Ric\wedge \id +(n-1)\epsilon \psi
\I\right)_{11}+\left(2\Ric\wedge \id +(n-1)\epsilon \psi
\I\right)_{22}\right).
\end{eqnarray*}
Here in the last equation above we have used Lemma 2.1 of \cite{BW},
which asserts that $R+R\#\I=\Ric\wedge \id$ (the use is not really
necessary). Since $\mu_1+\mu_2\ge 0$ and $\mu_\gamma \ge 0$ for all
$\gamma\ge 2$,
$$
\sum
(c^2_{1\beta\gamma}+c^2_{2\beta\gamma})\mu_\beta\mu_\gamma=2\sum_{\gamma
\ge
3}(c^2_{12\gamma}+c^2_{21\gamma})(\mu_1+\mu_2)\mu_\gamma+\sum_{\beta,
\gamma \ge
3}(c^2_{1\beta\gamma}+c^2_{2\beta\gamma})\mu_\beta\mu_\gamma\ge 0.
$$
Notice also that at $(x_0, t_0)$ we have that $\mu_{11}+\mu_{22}=0$,
which implies that $R_{11}+R_{22}=-2\epsilon \psi$, then $2\epsilon
\psi(x_0, t_0)\le 2A$. Hence at $(x_0, t_0)$ we have that
\begin{eqnarray}\label{str1}
0&\ge& \left(\frac{\p}{\p t} -\Delta \right)\left(\tilde
R_{11}+\tilde R_{22}\right)\nonumber\\
&\ge &\mu_1^2 +\mu_2^2+2\epsilon C_1\psi-200nA\epsilon \psi.
\end{eqnarray}
This is a contradiction if we choose $C_1>100nA$.
\end{proof}

By choosing the barrier function more carefully as in \cite{NT2,
N04} (see for example Theorem 2.1 of \cite{N04}), we can have the
following strong maximum principle.

\begin{corollary}\label{strong1}  Assume that $R(g(0))$ is
$2$-nonnegative and $2$-positive somewhere. Then there exists
$f(x,t)>0$ for $t>0$ and $f(x,0)\ge \frac{1}{2}(\mu_1+\mu_2)$, such
that $$\left(\mu_{1}+\mu_{2}\right)(x,t)\ge f(x, t).$$ In
particular, if $R(g(0))$ is $2$-nonnegative and $(\mu_1+\mu_2)(x_0,
t_0)=0$ for some $t_0\ge 0$, then $(\mu_1+\mu_2)(x, t)\equiv 0$ for
all $(x, t)$ with $t\le t_0$. Moreover, $\mu_1(x, t)=\mu_2(x,t)=0$
for all $(x, t)$ with $t\le t_0$ and
$$
\mathcal{N}_2(x, t)=\mbox{span}\{\omega_1, \omega_2\}
$$
is a distribution on $M$ which is invariant under the parallel
translation.
\end{corollary}

The above result together with (\ref{str1}) implies the following
classification of closed $2$-nonnegative manifolds.

\begin{corollary}\label{topo}  Assume that $R(g(0))$ is
$2$-nonnegative. Then for $t>0$, either the curvature operator
$R(g(t))$ is $2$-positive, or $R(g(t))\ge 0$. Hence Suppose $\left(
M^{n},g_{0}\right)  $ is a closed Riemannian manifold with
$2$-nonnegative curvature operator. Let $\tilde{g}\left(  t\right) $
be the lift to the universal cover $\tilde{M}$ of the solution
$g\left( t\right)  $ to the Ricci flow with $g\left(  0\right)
=g_{0}.$ Then for any $t>0$ we have either $\left(
\tilde{M}^{n},\tilde{g}\left( t\right)  \right)  $ is a closed
manifold with $2$-positive curvature operator or it is isometric to
the product of the following:

\begin{enumerate}
\item Euclidean space,

\item closed symmetric space,

\item closed Riemannian manifold with positive curvature
operator,

\item closed K\"{a}hler manifold with positive curvature operator on real
$\left(  1,1\right)  $-forms.
\end{enumerate}
\end{corollary}
\begin{proof} It follows from the above corollary and Hamilton's classification result on the
solutions with nonnegative curvature operator. See for example
\cite{CLN}, Theorem 7.34.
\end{proof}

Topologically, it is  now known, by \cite{BW}, that simply-connected
$2$-positive manifolds is  sphere, and the K\"ahler manifold in the
last case is biholomorphic to the complex projective space by the
earlier result of Mori -Siu-Yau. The fact that the curvature
operator of the evolving metrics becomes either $2$-positive or
nonnegative has been observed in \cite{Chen}. However, in
\cite{Chen} there is no clear statement of the strong maximum
principle, namely Corollary \ref{strong1}, on which the observation
relies. If evoking Theorem 2.3 of \cite{N04}, the splitting result
on solutions of Ricci flow on complete Riemannian manifold with
nonnegative curvature operator, we can write a similar statement
even when  $M$ is not assume to be compact. However, in this case
the Euclidean factor is only topological  (not isometric). Also we
do not know if a complete noncompact $2$-positive Riemannian
manifold is diffeomorphic to $\R^n$ or not.

\section{Manifolds with pinched curvature}

In \cite{H91} Hamilton proved that any convex hypersurface (with
dimension $\ge 3$) in Euclidean space with  second fundamental form
$h_{ij}\ge \delta \frac{\tr(h)}{n}\id$ must be compact.  In
\cite{CZ}, using the pre-established estimates of \cite{Hu} and
\cite{Sh2}, Chen and Zhu
 proved the following weak  version of above-mentioned Hamilton's
 result in terms of curvature operators. Namely, they
proved that {\it if a complete Riemannian manifold $(M^n, g)$ with
bounded and $(\epsilon, \delta_n)$- pinched curvature operator (with
$ n\ge 3$) in the sense that
$$
|R_{\W}|^2+|R_{\Ric_0}|^2\le \delta_n (1-\epsilon)^2|R_I|^2=\delta_n
(1-\epsilon)^2\frac{2}{n(n-1)}{\Scal (R)}^2
$$
for $\e>0$, $\delta_3>0, \delta_4=\frac{1}{5},
\delta_5=\frac{1}{10}$ and $\delta_n =\frac{2}{(n-2)(n+1)}$, where
$R_{\W}$, $R_{\Ric_0}$ and $R_{\I}$ denote the Weyl curvature
tensor, traceless Ricci part and the scalar curvature part. Then $M$
must be compact.} The strong pinching condition was the one
originally assumed in \cite{Hu}  to obtain various estimates and the
smooth convergence result. It was also shown in \cite{Hu} that it
implies that $R\ge \epsilon R_{I}$. In \cite{N05} the first author
showed that the above result of Chen-Zhu can be shown by the blow-up
analysis of \cite{H90} and some non-existence results on gradient
steady and expanding solitons obtained in \cite{N05}. (The detailed
proof on these non-existence results were submitted to 2004 ICCM
proceedings a while ago. See also forthcoming book \cite{CLN}.) With
the help of a family of invariant cones constructed in \cite{BW}, we
can now prove the following general result.

\begin{theorem}\label{ham}
Let $(M^n, g_0)$ be a complete Riemannian manifold with $n\ge 3$.
Assume that the curvature operator of $M$ is uniformly bounded
($|R_{ijkl}|(x)\le A$) and satisfies that
\begin{equation}\label{pin1}
R\ge \delta R_{\I}>0
\end{equation}
for some $\delta>0$. Then $(M, g)$ must be compact.
\end{theorem}

Recall that $R_{\I}=\frac{1}{n(n-1)}\Scal(R)\I$, where $\I$ is the
identity of $S_B^2(so(n))$. The above result is a natural  analogue
of Hamilton result for hypersurfaces.
\begin{proof} Let $g(t)$ be the solution to Ricci flow with initial metric $g_0$ constructed
by \cite{Sh1}. First we show that if $M$ is noncompact, $g(t)$ can
be extended to a long-time solution defined on $M\times [0,
\infty)$. In order to do that we first show that for sufficient
small $b>0$, $R(g_0)$ lies inside the invariant cone constructed by
Lemma 3.4 of \cite{BW}. Recall from \cite {BW} the linear
transformation
$$
\l_{a, b}: R\to R+2(n-1)aR_{\I}+(n-2)b R_{\Ric_0}.
$$
More precisely
\begin{eqnarray*}
\l_{a, b}(R)&=&R+2a \bar{\lambda}\I +2b\id \wedge \Ric_0(R)\\
&=& (1+2(n-1)a)R_{\I}+(1+(n-2)b)R_{\Ric_0}+R_{\W}.
\end{eqnarray*}
It is easy to see that $\l_{a,b}(S_B^2(so(n)))\subset S_B^2(so(n))$
and is invertible if $a\ne -\frac{1}{2(n-1)}$ and $b\ne
-\frac{1}{n-2}$. Using this linear map and Theorem 2 of \cite{BW}, a
pinching  family of invariant convex cones are constructed. In
particular, as one step of the construction, it was shown that
\begin{lemma}[B\"ohm-Wilking]\label{keylemma1}For $b\in [0, \frac{1}{2}]$, let
$$
a=\frac{(n-2)b^2+2b}{2+2(n-2)b^2}\, \mbox{ and }
p=\frac{(n-2)b^2}{1+(n-2)b^2}.
$$
Then the set $\l_{a, b}(C(b))$ where
$$
C(b)=\left\{R\in S_B^2(so(n))\, |\, R\ge 0,\, \Ric \ge
p(b)\frac{\tr(\Ric)}{n}\right\},
$$
is invariant under the vector fields $Q(R)$. In fact for $b\in(0,
\frac{1}{2}]$ it is transverse to the boundary of the set at all
boundary points $R\ne 0$.
\end{lemma}
We claim that there exists $b>0$ so small that $R(g_0)\in
\l_{a,b}(C(b))$, which is equivalent to  that $\l_{a,
b}^{-1}(R(g_0))\in C(b)$. For simplicity let $\tilde R=R(g_0)$,
$\bar{\lambda}(\tilde R)=\frac{\Scal(\tilde R)}{n}$ and $\l=\l_{a,
b}$. Direct computation shows that
$$
R:=\l^{-1}(\tilde
R)=\tilde{R}_{\W}+\frac{1}{1+2(n-1)a}\tilde{R}_{\I}+\frac{1}{1+(n-2)b}\tilde{R}_{\Ric_0}
$$
which implies that
$$
\Ric(\l^{-1}(\tilde R))=\frac{\bar{\lambda}(\tilde
{R})}{1+2(n-1)a}\id+\frac{1}{1+(n-2)b}\Ric_0(\tilde R)
$$
and
$$
\bar{\lambda}(R):=\frac{\tr(\l^{-1}(\tilde
R))}{n}=\bar{\lambda}(\tilde R)\left(1-\frac{2(n-1)
a}{1+2(n-1)a}\right).
$$
Let $\tl_i$ be the eigenvalues of $\Ric_0(\tilde R)$. Then by the
assumption (\ref{pin1}) we have that
\begin{equation}\label{pin2}
\tl_i+\bar{\lambda}(\tilde R)\ge \delta \bar{\lambda}(\tilde R).
\end{equation}
Clearly we also have that
\begin{equation}\label{pin3}
\tl_i+\bar{\lambda}(\tilde R)\le  n\bar{\lambda}(\tilde R).
\end{equation}
We first check that $R$ satisfies the Ricci pinching condition. In
fact if $\lambda_i$ are the eigenvalues of $\Ric_0(R)$, from the
above formulae we have that
\begin{eqnarray*}
-\lambda_i&=&-\frac{1}{1+(n-2)b}\tl_i\\
&\le& \frac{1-\delta}{1+(n-2)b}\bar{\lambda}(\tilde R)\\
&=& (1-\delta)\frac{1+2(n-1)a}{1+(n-2)b}\bar{\lambda}( R).
\end{eqnarray*}
Then there exist $\delta_1>0$ and  $b_0$ such that for all $b\in [0,
b_0]$, $-\lambda_i\le (1-\delta_1)\bar{\lambda}( R)$. Then we can
find $b_1\le b_0$ such that for any $b\in [0, b_1]$, $p(b)\le
\delta_1$. Hence  $R=\l_{a, b}^{-1}(\tilde R)$ satisfies the
pinching condition of $C(b)$. Now we check that $R=\l^{-1}_{a,
b}(\tilde R)\ge 0$. Rewrite
$$
R=\tilde{R}-\frac{2(n-1)a}{1+2(n-1)a}\tilde{R}_{\I}-\frac{(n-2)b}{1+(n-2)b}\tilde{R}_{\Ric_0}.
$$
Noticing that $a\to 0$ as $b\to 0$, we can find $b_2$ such that for
any $b\in [0, b_2]$ we have that
$$
R\ge \frac{\delta}{2}\tilde
R_{\I}-\frac{(n-2)b}{1+(n-2)b}\tilde{R}_{\Ric_0}.
$$
But the  eigenvalue (with respect to $e_i\wedge e_j$, where
$\{e_i\}$ is a basis of $TM$ consisting of eigenvectors of
$\Ric_0(\tilde R)$) of the right hand side operator can be computed
as
$$
\frac{\delta}{2}\frac{\bar{\lambda}(\tilde
R)}{n-1}-\frac{b}{1+(n-2)b}\left(\tl_i+\tl_j\right).
$$
Using (\ref{pin3}), the above can be bounded from below by
$$
\bar{\lambda}(\tilde
R)\left(\frac{\delta}{2(n-1)}-\frac{2(n-1)b}{1+(n-2)b}\right)>0
$$
if $b$ is close to $0$. This shows that there exists $b_3>0$ such
that for any $b\in (0, b_3]$, $R(g_0)\in \l_{a, b}(C(b))$.

Now the virtue of the proof of   Theorem 5.1 in \cite{BW}, along
with the short time existence result of \cite{Sh1}, shows that the
Ricci flow has long time solution. Otherwise, by Theorem 16.2 of
\cite{H90}, we would end up with a blow-up solution, which is
nonflat, noncompact, but whose curvature operator $R=R_{\I}$. In
view of Schur's theorem, this is a contradiction. Note that
$R(g_0)\in \l_{a, b}(C(b))$ allows us to apply the generalized
pinching set construction (Theorem 4.1)  from \cite{BW}, and since
the evolving metric has positive curvature operator and the manifold
is assumed to be noncompact, the injectivity radius always has a
lower bound in terms of the size of the curvature. All these
ingredients allow us to perform Hamilton's blow-up analysis
\cite{H90} (Theorem 16.2).

We continue to show that  the extra assumption that $M$ is
noncompact
 will lead us to a contradiction by performing the singularity analysis of
\cite{H90}  as $t\to \infty$. Notice that for all $t$, $R(g(t))$
will stay in the cone $\l_{a, b}(C(b))$ for some fixed (but
sufficiently small) $b$, by the tensor maximum principle, which can
be verified in the same way as Proposition \ref{chen1}. Now we claim
that
 the curvature of $g(t)$ satisfies that
\begin{equation}\label{ric-pin}
\Ric \ge p\frac{\tr(\Ric)}{n}\id.
\end{equation}
for some $p>0$. Let $R^*=R(g(t))$. First, by Lemma \ref{keylemma1}
we know that $R(g(t))\in \l_{a,b}(C(b)$ for some fixed small $b$.
Thus we can find $R\in C(b)$ such that $\l_{a, b}(R)=R^*$. Now let
$\bar{\lambda}=\frac{\tr(\Ric(R))}{n}$ and $\lambda_i$ be the
eigenvalues of $\Ric_0(R)$. By the assumption we have that
$-\lambda_i\le (1-p)\bar{\lambda}$. Now we compute the Ricci
curvature and its trace for $R^*$. By the definition of $\l_{a, b}$
we have that
$$
\Ric(R^*)=\Ric +2(n-1)a \bar{\lambda}\id +(n-2)b \Ric_0
$$
and
$$
\bar{\lambda}^*:=\frac{\tr(\Ric(R^*))}{n}=\bar{\lambda}(1+2(n-1)a).
$$
Letting $\lambda^*_i$ be the eigenvalue of $R^*$ we  have that $
\bar{\lambda}^*+\lambda_i^*
=(1+2(n-1)a)\bar{\lambda}+(1+(n-2)b)\lambda_i $. Therefore
\begin{eqnarray*}
-\lambda_i^*&=&-(1+(n-2)b)\lambda_i\\
&\le &(1-p)(1+(n-2)b)\bar{\lambda}\\
&=& (1-p)\frac{1+(n-2)b}{1+2(n-1)a}\bar{\lambda}^*\\
&\le &(1-p)\bar{\lambda}^*.
\end{eqnarray*}
Here we have used the fact that
$1+2(n-1)a=1+(n-1)\frac{(n-2)b^2+2b}{1+(n-2)b}>1+(n-2)b$. This
completes the proof of the claim (\ref{ric-pin}).

Since for all $g(t)$, its Ricci curvature satisfies (\ref{ric-pin}),
this holds up on the blow-down/blow-up solutions, which after
passing to its universal cover, are either a nonflat gradient steady
soliton or a nonflat gradient expanding soliton, with nonnegative
curvature operator, by results from \cite{H90} (Theorem 16.2,
Corollary 16.4) (See also  \cite{N02}, Theorem 4.2 and \cite{CZ}).
This contradicts to Corollary 3.1 of \cite{N05}.
\end{proof}

\section{Discussions}

In \cite{W}, the topology of so-called $p$-positive manifolds was
studied. In view of the result of B\"ohm-Wilking, it is reasonable
speculate that any noncompact complete Riemannian manifold with
$2$-positive curvature operator must be diffeomorphic to $\R^n$. In
\cite{N05} we speculated that any complete Riemannian manifolds with
positive pinched Ricci curvature must be compact. Theorem \ref{ham}
confirms it under stronger assumption on curvature operator.  The
problem in full generality  still remains unknown.

\medskip

{\it Acknowledgement}. Part of this paper was completed during the
first author's visit of ETH, Z\"urich. He would like to thank ETH,
especially Tom Ilamnen, for providing a stimulating environment and
various discussions. He also held informal discussions on \cite{BW}
with Ben Chow and Nolan Wallach.

\end{document}